\title{Algebraic Models for Homotopy Types}
\author{{\small \em Julio Rubio and Francis Sergeraert}}
\date{\footnotesize\today}

\documentclass[a4paper,12pt]{article}
\usepackage{amsfonts}
\usepackage{amssymb}

\newtheorem{thr}{Theorem}
\newtheorem{mth}[thr]{Meta-Theorem}

\newtheorem{dfn}[thr]{Definition}

\newtheorem{prb}[thr]{Problem}
\newtheorem{trm}[thr]{Terminology}
\newtheorem{rst}[thr]{Restriction}

\newcommand{\CCS}{\mathcal{CC}_{\mathcal{S}}}
\newcommand{\SSEH}{\mathcal{SS}_{EH}}

\newcommand{\bN}{\mathbb{N}}
\newcommand{\bR}{\mathbb{R}}
\newcommand{\bZ}{\mathbb{Z}}
\newcommand{\boxtt}[1]{\mbox{\small\texttt{#1}}}
\newcommand{\cA}{\mathcal{A}}
\newcommand{\cC}{\mathcal{C}}
\newcommand{\cH}{\mathcal{H}}
\newcommand{\cI}{\mathcal{I}}
\newcommand{\cO}{\mathcal{O}}
\newcommand{\cP}{\mathcal{P}}
\newcommand{\cS}{\mathcal{S}}

\newcommand{\cU}{\mathcal{U}}

\newcommand{\lrdc}{\mbox{\,\(\Leftarrow\hspace{-7pt}\Leftarrow\hspace{-7pt}\Leftarrow\)\,}}
\newcommand{\rrdc}{\mbox{\,\(\Rightarrow\hspace{-7pt}\Rightarrow\hspace{-7pt}\Rightarrow\)\,}}
\newcommand{\eqvl}{\mbox{\(\lrdc\hspace{-10pt}\rrdc\)}}

\newcommand{\bmp}{\rule{0pt}{1pt}\\[-15pt]{\tiny.\dotfill.}\\*[-3pt]}
\newcommand{\bmpi}{\begingroup\footnotesize}
\newcommand{\empi}{\endgroup\rule{0pt}{1pt}\\*[-2pt]}
\newcommand{\empim}{\endgroup\hspace{5pt}\mbox{\small\(\maltese\)}\rule{0pt}{1pt}\\*[-2pt]}
\newcommand{\empix}{\endgroup\rule{0pt}{1pt}\\[-2pt]}

\newcommand{\emp}{\rule{0pt}{1pt}\\*[-17pt]{\tiny.\dotfill.}\\*[-10pt]\rule{0pt}{1pt}}

\begin{document}

\voffset=-2.5cm
\hoffset=-0.9cm
\sloppy

\maketitle\vspace{-0pt}

\hfill\fbox{\parbox[t]{0.42\textwidth}{\footnotesize{\emph{
 \hspace*{0pt}\hfill As yet we are ignorant\\
 \hspace*{0pt}\hfill of an effective method of computing\\
 \hspace*{0pt}\hfill the cohomology of a Postnikov complex\\
 \hspace*{0pt}\hfill from
 \(\pi_n\) and \(k^{n+1}\)~\emph{\cite{EDM}}.
}}}}\vspace{20pt}

\begin{abstract}
The classical problem of \emph{algebraic models} for homotopy types is
\emph{precisely stated}, to our knowledge for the first time. Two different
natural statements for this problem are produced, the simplest one being
\emph{entirely solved} by the notion of \(\SSEH\)-structure, due to the
authors. Other tentative solutions, Postnikov towers and \(E_\infty\)-chain
complexes, are considered and compared with the \(\SSEH\)-structures. In
particular, which looks at least like an unfortunate imprecision in the usual
definition of the \(k\)-``invariants'' is explained, which implies we seem far
from a solution for the ideal statement of our problem. At the positive side,
the problem stated\footnote{Probably badly translated from Japanese; it must be
understood at the last line: ``from~\((\pi_n, k^{n+1})_{n \geq 2}\)''.} above
in the framed quotation is solved.
\end{abstract}

\section{Introduction\protect\footnote{This
text is a slightly expanded version of a talk given by the authors at the
RSME-AMS meeting at Sevilla in June 2003, which explains its nature a little
expository. The talk itself is available as a dvi-file at \boxtt{\footnotesize
www-fourier.ujf-grenoble.fr/\~{ }sergerar/Talks}.}.}

Obtaining ``algebraic'' models for \(\mathbb{Z}\)-homotopy types is a major
problem. The \emph{statement} of the problem itself is a constant source of
strong and regrettable ambiguities. We explain in this article why the
adjective \emph{algebraic} is in fact inappropriate, the right one being
\emph{computable} (or effective, constructive, \ldots).

The problem of the title can \emph{then} be precisely stated in two different
ways, the \emph{hard} problem (Problem~\ref{vkkol} in Section~\ref{aajnx}) and
the \emph{soft} problem (Problem~\ref{srxvi} in Section~\ref{ydxne}). The
notion of simplicial set with effective homology~(\(\SSEH\)), due to the
authors, is a complete solution for the soft problem, very simple from a
theoretical point of view, once the possibilities of functional programming are
understood. This solution has led to an interesting concrete computer work, the
Kenzo program, demonstrated a little in the article to give to the reader an
\emph{experimental} evidence that the stated results are correct.

Other solutions for the soft problem could be based on operadic techniques, and
they are now intensively looked for. The key point is the notion of
\(E_\infty\)-operad; a broad outline of the main results so obtained is given
and compared with the \(\SSEH\) solution. The current result is that the
\(\SSEH\) solution is, for the soft problem, terribly simple; furthermore it
is\ldots\ available. The operadic structures are of course interesting, give
many useful informations, but are \emph{by-products} of \(\SSEH\)-structures;
furthermore it is not clear how they could produce \emph{autonomous} computable
objects. The good point of view for future work is probably a mixture of
\(\SSEH\)'s and operadic techniques, the last ones to be considered as good
tools to better understand and also to improve the computability results so
easily obtained through \(\SSEH\)'s.

The \emph{hard problem} is \emph{so} reduced to the problem of equivalence
between sets of \(k\)-invariants, problem which, up to further information,
seems open: we explain why the so-called \(k\)-invariants are not actual
invariants so that finally the standard Postnikov theory \emph{does not} solve
the \emph{hard problem}.

\section{The right statement of the problem.}\label{aajnx}

The construction of \emph{algebraic models} for \emph{homotopy types} is a
``classical'' problem in Algebraic Topology which, to our knowledge, has never
been precisely stated, that is, \emph{mathematically} stated. Experience shows
the topologists have a rather imprecise idea about the exact nature of this
problem, a situation frequently leading to misunderstandings or even sometimes
to severe imprecisions; an example of this sort being the usual belief that the
so-called \(k\)-invariants are\ldots\ invariants, an erroneous appreciation,
see Sections~\ref{mssfy} and~\ref{deoad}.

Most of the topologists should agree with the following statement of our
problem.

\begin{prb}\label{uhvfr} ---
Let \(\cH\) be the \emph{homotopy category}. How to design an \emph{algebraic
category} \(\cA\) and a functor \(F: \cH \rightarrow \cA\) which is an
equivalence of categories?
\end{prb}

Instead of working in the category \(\cH\), reputed to be a difficult category,
you might work in the category \(\cA\), an algebraic category, hence probably a
more convenient workspace. The image \(F(X)\) of some homotopy type \(X\) would
be an \emph{algebraic} object, for example a chain complex provided with a
sufficiently rich structure to entirely define a homotopy type.
Problem~\ref{uhvfr} leads to an auxiliary problem.

\begin{prb}\label{atzir} ---
What is the definition of an \emph{algebraic} category?
\end{prb}

It happens that standard logic shows such a definition \emph{cannot exist};
this is a direct consequence of the formalization of mathematics, asked for by
Hilbert, and realized through various systems, mainly the so-called
Zermelo-Fraenkel and Bernays-Von Neumann systems. In a sense, formalization of
mathematics consists in making \emph{entirely} algebraic our mathematical
environment, even when we work in fields that are not usually considered as
algebraic, like in analysis, probability, and also in topology.

The following example is fairly striking. Most of the topologists consider a
simplicial set is \emph{not} an algebraic object. A simplicial set \(S\) is a
sequence of simplex sets \((S_n)\) combined with some sets of operators between
these simplex sets, appropriate composites of these operators having to satisfy
a few simple relations. Most of the topologists consider a chain complex
\(C_\ast\) provided with a module structure with respect to some (\ldots
algebraic!) operad \(\cO\) \emph{is} an algebraic object. Such a chain complex
is a sequence of chain groups \((C_n)\) combined with some sets of operators
between these chain groups and their tensor products, appropriate composites of
these operators having to satisfy a large set of sophisticated relations. Where
is the basic difference? This appreciation --- an \(\cO\)-module \emph{is} an
algebraic object and a simplicial set \emph{is not} --- is arbitrary.
Furthermore a simplicial structure is simpler than an \(\cO\)-module structure,
so that a beginner in the subject would probably guess the first structure type
is ``more'' algebraic than the second one. Must we recall we are working in
mathematics, not in philosophy? Our workspace require \emph{mathematical}
definitions, not fuzzy speculative claims based only on vague traditions.

\begin{trm}\label{aytus} ---
In our current mathematical environment, the border between \emph{algebraic}
objects and \emph{non-algebraic} objects cannot be \emph{mathematically}
defined.
\end{trm}

Let us continue our comparison between simplicial sets and chain complexes,
which will eventually lead to the right point of view. The simplest example of
an interesting result produced by Algebraic Topology is the Brouwer theorem, a
direct consequence of the following.

\begin{thr} ---
Let \(i_n: S^{n-1} \rightarrow D^n\) be the canonical inclusion of the
\((n-1)\)-sphere into the \(n\)-ball. There does not exist a continuous map
\(\rho_n: D^n \rightarrow S^{n-1}\) such that the composite \(\rho_n \circ
i_n\) is the identity map of \(S^{n-1}\).
\end{thr}

In fact, if you apply the \(H_{n-1}\)-functor to the data, the statement is
transformed into: let \(i: \bZ \rightarrow 0\) be the null morphism; there does
not exist a morphism \(\rho: 0 \rightarrow \bZ\) such that the composite \(\rho
\circ i\) is the identity morphism of \(\bZ\).

Most of the topologists think this process produces the result because the
transformed problem has an \emph{algebraic} nature, but this is erroneous. The
\emph{algebraic} qualifier is secondary and, as previously explained, cannot be
mathematically justified. The right qualifier in fact is \emph{computable}. The
transformed problem is a particular case of the following: let \(m\), \(n\) and
\(p\) be three non-negative integers, and \(f: \bZ^m \rightarrow \bZ^n\) and
\(F: \bZ^m \rightarrow \bZ^p\) be two \(\bZ\)-linear morphisms; does there
exist a morphism \(g: \bZ^n \rightarrow \bZ^p\) satisfying \(g \circ f = F\)?
It is common to think of this problem as an algebraic one, but in fact the only
important point for us is that there exists an \emph{algorithm} giving the
solution: a Smith reduction of the \(\bZ\)-matrices representing \(f\) and
\(F\) quickly gives the solution; in the case of the Brouwer problem, the Smith
reduction is already done.

The previous considerations about simplicial sets give another idea. Because a
simplicial set is in fact as ``algebraic'' as a homology group or a chain
complex, why not work directly with simplicial models for \(S^{n-1}\) and
\(D^n\)? It is easy to give simplicial models with two (resp. three)
non-degenerate simplices for \(S^{n-1}\) (resp. \(D^n\)), models that are
undoubtedly ``algebraic''. But these models have an essential failing: they do
not satisfy the Kan extension condition, so that they are not appropriate for
working in the homotopy category \(\cH\).  In general the Kan simplicial models
are highly infinite and cannot be directly used for computations: any tentative
solution using in an essential way the Kan simplicial sets raises hard
computability problems. We will see later that our solution for ``algebraic''
models for homotopy types is a simple but subtle combination of simplicial sets
most often \emph{not of finite type} with chain complexes of finite type.

There is a common fundamental confusion between the \emph{algebraic} and
\emph{computable} qualifiers, still present in the ordinary understanding of
the very nature of Algebraic Topology. From this point of view, it can be
useful to recall the frequent opinion of the pupils in secondary schools: ``I
prefer Algebra rather than Geometry, because in Algebra we can use
\emph{automatic} methods giving the results that are looked for; on the
contrary, in Geometry, we often have to \emph{discover} the appropriate method
for some particular problem''; another example of the same confusion between
\emph{algebraic} and \emph{computable}.

Let us look again at the statement of Problem~\ref{uhvfr}. We see the
requirement for the category~\(\cA\) to be algebraic cannot be defined; in fact
we are looking for a target category where \emph{automatic} computations
(pleonasm) can be undertaken. We so obtain a new statement for our problem.

\begin{prb}\label{vkkol} \textbf{\emph{(Hard Problem)}} ---
Let \(\cH\) be the \emph{homotopy category}. How to design a \emph{computable
category} \(\cC\) and a functor \(F: \cH \rightarrow \cC\) which is an
equivalence of categories?
\end{prb}

With the satellite problem:

\begin{prb} ---
What is the definition of a \emph{computable} category?
\end{prb}

We do not want to consider the details of the last subject, an interesting
subject, out of scope of the present paper; fundamentally different answers are
possible, mainly from the following point of view: do you intend to apply the
``computable'' qualifier to the \emph{elements} of an object in the category or
to the \emph{objects} themselves, or both? To our knowledge, the relevant
corresponding theory is not yet settled\footnote{For example the
reference~\cite{RDBR}, interesting, cannot be useful for our main problem; look
for the entry \emph{equality} in the index, and you will quickly understand
that no tool is provided there for the equality problem between objects of a
category.}. The few examples given in the paper could be a guideline toward the
most natural solutions of this question.

In this paper, a \emph{computable category} is a category having properties
roughly similar to those that are exhibited for our \(\SSEH\) category, and
this approximate ``definition'' is here sufficient.

In other words, a click on the {\small \texttt{rename}} button, replacing the
imprecise identifier \emph{Algebraic Topology} by the precise one
\emph{Constructive Topology}, could be a good idea.

\section{Three tentative solutions.}\label{ydxne}

The current state of Algebraic Topology gives mainly three possibilities:
\begin{enumerate}
\item
The Postnikov category;
\item
The authors' solution: the category \(\SSEH\);
\item
The operadic solutions.
\end{enumerate}

In short, the first possibility is currently inadequate in the standard
framework, because of an essential lack of computability, see the framed title
inscription, and also because of the underlying classification problem which
does not yet seem solved. The \(\SSEH\) category solves a \emph{subproblem},
the ``soft'' problem stated a little later, and furthermore makes the Postnikov
category computable; a consequence is the fact that the Postnikov category,
when modelled as a satellite category of the \(\SSEH\) category, solves the
same subproblem. It can be reasonably conjectured that the third idea, using
operadic techniques, should in finite time solve the same subproblem, but we
are still far from it, and this seems the challenge \#1 for the operad
developers: how to organise the \(E_\infty\)-chain complexes as an autonomous
computable category? The impressive \emph{concrete} results obtained by V.
Smirnov~\cite{SMRN}, at least when working with coefficients in a field
\(\mathbb{F}_p\), indicate that the vast theoretical study about
\(E_\infty\)-chain complexes undertaken by this author could be a good
guideline.

 The gap about the classification problem remains present for
the three techniques.

Once the theoretical and concrete possibilities of functional programming are
understood, the \(\SSEH\) category is not complicated, so that it has been
possible to write down a computer program implementing the \(\SSEH\) category
and to use it, see~\cite{DRSS} and Sections~\ref{czlgy} and~\ref{jmahz} of the
present paper.

\subsection{The Postnikov category.}\label{mssfy}\vspace{15pt}

\begin{rst} ---
Unless otherwise stated, all our topological spaces are connected and
simply connected.
\end{rst}

An object of the Postnikov category is a pair of sequences \(((\pi_n)_{n \geq
2}, (k_n)_{n \geq 3})\), made of homotopy groups and ``\(k\)-invariants''
defining a \emph{Postnikov tower} \((X_n)_{n \geq 2}\). The first stage of the
tower \(X_2\) is \(K(\pi_2, 2)\), the first \(k\)-invariant \(k_3 \in H^4(X_2,
\pi_3)\) defines a fibration \(X_3 \twoheadrightarrow X_2\) the fiber of which
being \(K(\pi_3, 3)\), and so on. It is not hard to define the valid morphisms
between two towers, and we have so defined the Postnikov category \(\cP\). We
know it is not a common opinion, but this category is as ``algebraic'' as the
usual so-called algebraic categories; notice in particular the ingredients
defining a Postnikov tower are commutative groups and elements of some
commutative groups; are not they \emph{algebraic}?

The so-called \(k\)-invariants are not invariants, for the following reason:
\emph{different} \(k\)-invariants frequently give the \emph{same} homotopy
type. Identifying the corresponding equivalence classes is a problem which, to
our knowledge, is yet without any solution. Let us look at this simple example:
what about the Postnikov towers with only \(\pi_2 = \bZ^p\), \(\pi_5 = \bZ\)
and the other \(\pi_n\)'s are null. The only relevant \(k\)-invariant is \(k_5
\in H^6(K(\pi_2, 2), \pi_5) = \mbox{Cub}(\bZ^p, \bZ)\), the \(\bZ\)-module of
the cubical forms over \(\bZ^p\); making these cubical forms actual invariants
amounts to being able to construct and describe in a computational way the
quotient set \(\mbox{Cub}(\bZ^p, \bZ)\,/\,\mbox{(linear equivalence)}\). We
have questioned several arithmeticians and they did not know whether
appropriate references would allow a \(k\)-invariant user to solve this
problem: the classification problem does not seem to be \emph{solved} by the
``\(k\)-invariants'' and our example is one of the simplest ones\footnote{We
would like to thank Daniel Lazard for his study (private communication) which
opens several interesting research directions around this subject.}.

Let us quote certainly one of the best specialists in homotopy theory. Hans
Baues explains in \cite[p. 33]{BAUS4}:

{\newlength{\parquot}
 \setlength{\parquot}{\textwidth}
 \addtolength{\parquot}{-\parindent}
 \addtolength{\parquot}{-\parindent}
 \addtolength{\parquot}{-\parindent}
 \begin{center}
 \fbox{\parbox{\parquot}
 {Here \(k_n\) is actually an \emph{invariant} of the homotopy type of \(X\) in
 the sense that a map \(f: X \rightarrow Y\) satisfies
 \[
 (P_{n-1}f)^\ast k_n Y = (\pi_n f)_\ast k_n X
 \]
 in \(H^{n+1}(P_{n-1} X, \pi_n Y)\).}}
 \end{center}}

This explanation is not correct; the cohomology class \(k_n\) would be an
actual \emph{invariant} of the homotopy type if a homotopy equivalence
\mbox{\(f: X \simeq Y\)} implies \(k_n X\ \fbox{=}\ k_n Y\); in fact the framed
equal sign does not make sense: the underlying cohomology groups are not the
\emph{same}, they are only, in the relevant cases, \emph{isomorphic} and two
invariants should be considered as ``equal'' as soon as they are in turn
``isomorphic'' in an obvious sense. Baues' relation only shows the
``\(k\)-invariant'' depends \emph{functorially} on the data, but it is not an
\emph{invariant}; the definition would be acceptable if the isomorphism problem
between the various possible \(k\)-invariants in the same homotopy class had a
(computable) solution, but the simple example given before shows such a
solution does not seem currently known. We will examine again this question in
a more explicit way in Section~\ref{deoad}, where another classical reference
about \(k\)-invariants~\cite{MAY} is also studied.

This is probably the reason why Hans Baues uses entirely different techniques
to obtain certainly the most interesting concrete results so far reached in the
\emph{general} classification problem; see~\cite[Section~11]{BAUS4} and Baues'
references in the same paper.

Let us consider the following subproblem of the hard one:

\begin{prb}\label{srxvi} \textbf{\emph{(Soft Problem)}} ---
How to design a computable category \(\cC\) and a functor \(F: \cC \rightarrow
\cH\) such that any recursive homotopy type is in the image of \(F\)?
\end{prb}

In fact the \emph{hard problem} as it is stated in Problem~\ref{vkkol} cannot
have a solution: the standard homotopy category \(\cH\) is much too rich to
make it equivalent to a computable category. This is a situation analogous to
which is well known for example for the real numbers. A computable real number
is usually called a \emph{recursive} real number and the set of the recursive
real numbers is countable, much smaller than the set of ``ordinary'' real
numbers, see~\cite{TRVN}. In the same way:

\begin{dfn}\label{gqifn} \emph{---
A \emph{recursive homotopy type} is defined by a recursive Postnikow tower
\(((\pi_n)_{n \geq 2}, (k_n)_{n \geq 3})\): the data of this tower are defined
by an \emph{algorithm} \(n \mapsto (\pi_n, k_n)\). In other words the recursive
homotopy category is the image of the canonical functor \(\cP_r \rightarrow
\cH\) if \(\cP_r\) is the category of the recursive Postnikov towers.}
\end{dfn}

In fact, in the standard context, this definition does not make sense: the
required algorithm must be able to \emph{compute} the \(H^{n+1}(X_{n-1},
\pi_n)\) to allow it to ``choose'' the next \(k_n\), and classical Algebraic
Topology does not solve this question~\cite{EDM}. To our knowledge, there are
currently only two solutions for this problem, independantly and simultaneously
found by Rolf Sch\"on~\cite{SCHN} and the present authors
\cite{SRGR3,RBSR6,SRGR5}. We will see later that our \(\SSEH\) category allows
us in particular to compute the cohomology groups \(H^{n+1}(X_{n-1}, \pi_n)\)
when the previous data are available, so recursively defining where the \(k_n\)
is to be chosen. In this way our category \(\SSEH\) makes coherent
Definition~\ref{gqifn} and, \emph{then only}, the Postnikov category becomes an
obvious solution for the \emph{soft problem}. In fact we will also see the
\(\SSEH\) category directly gives a solution for the \emph{soft problem}.

It should be clear now that in the statement of the \emph{hard problem}, the
category~\(\cH\) must obviously be replaced by the category of the
\emph{recursive} homotopy types. Up to a finite dimension, this amounts only to
requiring that the homotopy groups \(\pi_n\) are of finite type, but if the
situation is considered without any dimension limit, the requirement is much
stronger.

\begin{rst}\label{gmsev} ---
From now on, all our categories are implicitely limited to \emph{recursive}
objects and \emph{recursive} morphisms.
\end{rst}

The \emph{soft problem} then is the same as the \emph{hard problem} except that
the classification question is given up.

\subsection{The operadic solutions.}\label{bdpci}

Many interesting works have been and are currently undertaken to reach operadic
solutions for the \emph{hard} and for the \emph{soft problem}. Probably the
most advanced one is due to Michael Mandell~\cite{MNDL}, usually considered as
a ``terminal'' solution. Of course we do not intend to reduce the interest of
this work, essential, but Mandell's article \emph{is not} a solution for the
\emph{hard problem}: the computability question is not considered, and the
proposed solution invokes numerous layers of sophisticated techniques, so that
the computational gap is not secondary. In fact the interesting next problem
raised by Mandell's paper is the following: is it possible to ``naturally''
extend the paper to obtain the corresponding \emph{effective} statements, or on
the contrary is it necessary to add something which is essentially new and/or
different? A solution is probably reachable in characteristic \(p\), but if it
was possible to obtain the same result with respect to the ground ring \(\bZ\),
then crucial computability problems in arithmetic would be solved, problems
which cannot be directly reduced to \mbox{\(\bZ_p\)-problems}, see the
discussion in Section~\ref{nnmrf}.

The operadic techniques raise other essential difficulties: it seems extremely
difficult to make \emph{computable} the relevant categories. The challenge is
the following: the ordinary constructions of algebraic topology --- loop
spaces, classifying spaces, fibrations for example --- should have a
\emph{translation} in the chosen category. The now standard methods of
\emph{closed model categories} give many possibilities, give frequently elegant
theoretical solutions for these translations, but the computational satellite
problems are seldom studied, why? When we observe the terrible problems met by
the ``classical'' topologists when they try to iterate the cobar construction
in a purely algebraic framework, we cannot be very optimist. However, as
already observed, an entirely combinatorial translation of Smirnov's operadic
techniques could be the right direction. Another solution could consist in
\emph{combining} \(\SSEH\)-structures and operadic structures; this question
will be examined in Section~\ref{bxwta}.

\section{The category \(\SSEH\).}

\begin{rst} ---
All the chain complexes considered from now on are implicitly assumed to be
\emph{free} \(\bZ\)-complexes, not necessarily of finite type.
\end{rst}

The notion of \emph{reduction}\footnote{Often called \emph{contraction}, but
this is a non-negligible terminological error: a contraction is a
\emph{topological} object and a reduction is only an \emph{algebraic} object;
it is important to understand a reduction \emph{does not} solve the underlying
\emph{topological} problem. Exercise: why this remark does not contradict
Terminology~\ref{aytus}?} is well known.

\begin{dfn} \emph{---
A \emph{reduction} \(\rho: C_\ast \rrdc D_\ast\) between two chain complexes
\(C_\ast\) and~\(D_\ast\) is a triple \(\rho = (f,g,h)\) where:
\begin{enumerate}
\item
The first component \(f\) is a chain complex morphism \(f: C_\ast \rightarrow
D_\ast\);
\item
The second component \(g\) is a chain complex morphism \(g: D_\ast \rightarrow
C_\ast\);
\item
The third component \(h\) is a homotopy operator (degree = +1) \(h: C_\ast
\rightarrow C_\ast\);
\item These components satisfy the relations:
\begin{enumerate}
\item
\(f \circ h = 0\);
\item
\(h \circ g = 0\);
\item
\(h \circ h = 0\);
\item
\(\mbox{id}_{D_\ast} = f \circ g\)
\item
\(\mbox{id}_{C_\ast} = g\circ f + d_{C_\ast} \circ h + h \circ d_{C_\ast}\).
\end{enumerate}
\end{enumerate}}
\end{dfn}

These relations express in an \emph{effective} way how the ``big'' chain
complex \(C_\ast\) is the direct sum of the ``small'' one \(D_\ast\) and an
\emph{acyclic} one, namely the kernel of \(f\).

\begin{dfn}\label{jvcgb} \emph{---
A \emph{strong chain equivalence} (or simply an \emph{equivalence}):
\[
 \varepsilon: C_\ast \eqvl D_\ast
\]
is a pair of reductions \(\varepsilon = (\rho_\ell, \rho_r)\) where:
\[
C_\ast \stackrel{\raisebox{5pt}{\(\rho_\ell\)}}{\lrdc} \widehat{C}_\ast
\stackrel{\raisebox{5pt}{\(\rho_r\)}}{\rrdc} D_\ast
\]
with \(\widehat{C}_\ast\) some intermediate chain complex.}
\end{dfn}

\begin{dfn} \emph{---
A \emph{simplicial set with effective homology} is a 4-tuple \[X_{EH} = (X,
C_\ast X, EC^X_\ast, \varepsilon^X)\] where:
\begin{enumerate}
\item
The first component \(X\) is a \emph{locally effective} simplicial set;
\item
The second component \(C_\ast X\) is the \emph{locally effective} chain complex
canonically associated to \(X\);
\item
The third component \(EC^X_\ast\) is an \emph{effective} chain complex;
\item
The last component \(\varepsilon^X\) is a \emph{strong chain equivalence}
\(\varepsilon^X: C_\ast X \eqvl EC^X_\ast\).
\end{enumerate}}
\end{dfn}

An \emph{effective} chain complex is an ordinary object, no surprise; it is an
algorithm \(n \mapsto (C_n, d_n)\) where, for every integer~\(n\), the
corresponding chain group \(C_n\) is a free \(\bZ\)-module of finite type, and
\(d_n\) is a \(\bZ\)-matrix describing the boundary operator \(d_n: C_n
\rightarrow C_{n-1}\). Elementary algorithms then allow to compute the homology
groups of such a complex. The third component~\(EC_\ast^X\) of a simplicial set
with effective homology is of this sort.

A \emph{locally effective} chain complex is \emph{quite different}. It is an
algorithm:
\[n \mapsto (\chi_n, d_n)\] to be interpreted as follows.
\begin{enumerate}
\item
The first component \(\chi_n\) of a result is also an algorithm \(\chi_n: \cU
\rightarrow \{\top, \bot\}\) where~\(\cU\) (universe) is the set of \emph{all}
the machine objects, so that for every machine object \(\omega\), the algorithm
\(\chi_n\) returns \(\chi_n(\omega) \in \{\top, \bot\}\), that is, true or
false, true if and only if \(\omega\) is a \emph{generator} of the \(n\)-th
chain group of the underlying chain complex.
\item
The second component \(d_n\) of a result is again an algorithm: if
\(\chi_n(\omega) = \top\), then \(d_n(\omega)\) is defined and is the boundary
of the generator \(\omega\), therefore a finite \(\bZ\)-combination of
generators of degree \(n-1\).
\end{enumerate}

The set \(\cU\), for any reasonable machine model, is infinite countable, so
that a locally effective chain complex in general \emph{is not of finite type}.
The adverb \emph{locally} has the following meaning: if someone produces some
(every!) generator \(\omega\) of degree~\(n\), then the \(d_n\)-component is
able to compute the boundary \(d_n(\omega)\). The terminology
\emph{generator-wise} effective chain complex would be more precise but a
little unwieldy.

A non-interesting but typical example of locally effective chain complex would
be produced by \(\chi_n(\omega) = \top\) if and only if \(\omega \in \bN\),
independently of \(n\), and \(d_n(\omega) = 0\) for every \(n \in \bZ\) and
\(\omega \in \bN\). In other words the underlying chain complex would be the
periodic one \(C_n = \bZ^{(\bN)}\) with a null boundary.

Standard logic shows in general the homology groups of a locally effective
chain complex \emph{are not} computable; this is an avatar of the
G\"odel-Turing-Church-Post theorems about incompleteness. More generally, a
\emph{global} information in general cannot be deduced from a locally effective
object. The second component \(C_\ast X\) of a simplicial set with effective
homology is of this sort.

A locally effective simplicial set is defined in the same way; the simplices
are defined through characteristic algorithms \(\chi_n\), and instead of
computing boundaries, a set of appropriate operators compute faces and
degeneracies.

The second component \(C_\ast X\) of a simplicial set with effective homology
is redundant: a simple algorithm can construct it from the locally effective
simplicial set \(X\); and strictly speaking, we could forget it in the
presentation. But the key points in an object with effective homology are:
\begin{enumerate}
\item
The main components are two \(\bZ\)-free chain complexes \(C_\ast X\) and
\(EC_\ast^X\), the first one being a direct consequence of the underlying
object \(X\), the second one describing the homology of this object, reachable
through an elementary algorithm;
\item
The component \(C_\ast X\) is \emph{locally effective} allowing it not to be of
finite type, with the drawback that in general its homology is not computable;
\item
The component \(EC_\ast^X\) is \emph{effective}, therefore of finite type with
a computable homology;
\item
The equivalence \(\varepsilon^X\) is the key connection between the locally
effective object \(C_\ast X\) and the effective one \(EC_\ast^X\).
\end{enumerate}
and it is hoped the nature of this organization is better explained in the
notation \((X, C_\ast X, EC_\ast^X, \varepsilon^X)\).

Let us insist again on a key point: the notation is very misleading; the
\emph{locally effective} subobject \(X\) which is the first component of an
object with effective homology \emph{does not \underline{effectively}}
determine the ``mathematical'' underlying object \(X\), because of the standard
incompleteness theorems. Only very partial --- ``local'' --- informations are
reachable through such an object; if free colors were available in the
text-processing system used when preparing this text, a very pale color should
have been chosen for this symbol \(X\), to clearly recall this sub-object
\emph{is not} \(X\), but a new kind of object rarely considered in standard
mathematics, an object \emph{of the third type}~\cite{SRGR5}.

\begin{thr}\label{vntdw} ---
The category \(\SSEH\) is a solution for the \emph{soft problem}.
\end{thr}

It is not possible in the framework of this paper to give a \emph{proof} of
Theorem~\ref{vntdw}, we will give only a \emph{demonstration}. We apologize for
the poor joke: ``demonstration'' has two different meanings in our context, it
can be a \emph{mathematical} proof, and it can be also a \emph{machine}
(computer) demonstration. It is expected in this case a machine demonstration
should give to the reader a strong conviction the Kenzo program \emph{contains}
a proof of Theorem~\ref{vntdw}. This is the aim of Sections~\ref{czlgy}
and~\ref{jmahz}.

\section{A small machine demonstration.}\label{czlgy}

This section uses a small machine demonstration to explain how, thanks to the
powerful computer language Common Lisp, the Kenzo program\cite{DRSS} makes the
objects and morphisms of the \(\SSEH\) category concretely available to the
topologist.

Let us consider the following space:
\[{X} = {\Omega(\Omega (P^\infty(\bR)/P^3(\bR)) \cup_4 D^4) \cup_2
 D^3}\]

The infinite real projective space truncated to the dimension 4,
\(P^\infty(\bR)/P^3(\bR)\), is firstly considered; its loop space is
constructed and the homotopy of this loop space begins with \(\pi_3 = \bZ\); so
that attaching a 4-cell by a map \(\partial D^4 \rightarrow S^3\) of degree~4
makes sense and this is done. The loop space functor is again applied to the
last space and finally a 3-cell is attached by a map of degree~2. This
artificial space~\(X\) is chosen because it is not too complicated, yet it
accumulates the main known obstacles to the theoretical and concrete
computation of homology groups in small dimensions.

The space \(X\) is an object of the category \(\SSEH\), so that the Kenzo
program can construct it as such an object. As follows:

 \bmp
 \bmpi\verb|> (progn|\empix
 \bmpi\verb|   (setf P4 (r-proj-space 4))|\empi
 \bmpi\verb|   (setf OP4 (loop-space P4))|\empi
 \bmpi\verb|   (setf attach-4-4|\empi
 \bmpi\verb|     (list (loop3 0 4 4) (loop3) (loop3) (loop3) (loop3)))|\empix
 \bmpi\verb|   (setf DOP4 (disk-pasting OP4 4 'D4 attach-4-4))|\empix
 \bmpi\verb|   (setf ODOP4 (loop-space DOP4))|\empix
 \bmpi\verb|   (setf attach-3-2|\empi
 \bmpi\verb|     (list (loop3 0 (loop3 0 4 1) 2) (loop3) (loop3) (loop3)))|\empi
 \bmpi\verb|   (setf X (disk-pasting ODOP4 3 'D3 attach-3-2)))|\empim
 \emp

We cannot explain the technical details of the construction, but most of the
statements are self-explanatory. Each object is located by a symbol and the
assignment is set through a \boxtt{setf} Lisp statement. For example the
initial truncated projective space is assigned to the symbol \boxtt{P4}. An
object such as \boxtt{attach-4-4} describes an attaching map as a
\emph{simplicial} map \(\partial \Delta^4 \rightarrow \boxtt{OP4}\) and this
description is then used by the Lisp function \boxtt{disk-pasting} which
constructs the desired space by attaching a cell according to the descriptor
\boxtt{attach-4-4}. The same for the end of the construction.

When this statement is executed, Lisp returns:

 \bmp
 \bmpi\verb|   ...|\empi
 \bmpi\verb|   (setf X (disk-pasting ODOP4 3 'D3 attach-3-2)))|\empim
 \bmpi\verb|[K17 Simplicial-Set]|\empi
 \emp

A maltese cross {\small\(\maltese\)} means the Lisp statement is complete, the
\boxtt{read} stage of the \boxtt{read-eval-print} Lisp cycle is finished, the
\boxtt{eval} stage starts for an execution of the just read Lisp statement, it
is the stage where the machine actually works, \emph{evaluating} the statement;
most often, an object is \emph{returned} (printed), it is the result of the
evaluation process, in this case the simplicial set \#\boxtt{K17}, located
through the~\boxtt{X} symbol. This object \boxtt{X} is a (machine) version
\emph{with effective homology} of the topological space \(X\).

So that we can ask for the \emph{effective homology} of \(X\); it is reached by
the function \boxtt{efhm} (\underline{ef}fective
\underline{h}o\underline{m}ology) and assigned to the symbol \boxtt{SCE}
(\underline{s}trong \underline{c}hain \underline{e}quivalence):

 \bmp
 \bmpi\verb|> (setf SCE (efhm X))|\empim
 \bmpi\verb|[K268 Equivalence K17 <= K256 => K258]|\empi
 \emp

The Kenzo program \emph{returns} a (strong) equivalence
(Definition~\ref{jvcgb}) between the chain complexes \#\boxtt{K17} and
\#\boxtt{K258}. Usually it is understood a simplicial set \emph{produces} an
associated chain complex, but we may conversely consider that a simplicial set
is nothing but a chain complex where a simplicial structure is \emph{added},
compatible with the differential; it is the right point of view and Kenzo
follows this idea. Please compare with the discussion after
Problem~\ref{atzir}: it should be more and more obvious that a simplicial set
is itself a chain complex with a further \emph{algebraic}~(!) structure. In
other words if you are only looking for an \emph{algebraic model for a homotopy
type}, the notion of simplicial set is a simple undeniable definitive solution,
already given fifty years ago by Eilenberg and MacLane~\cite{ELMC1,ELMC2}; this
is the reason why you \emph{must} add a computability requirement to finally
obtain an interesting problem.

In our equivalence describing the effective homology of \(X\), the right chain
complex \#\boxtt{K258} is \emph{effective}, the left one \#\boxtt{K17} is not:

 \bmp
 \bmpi\verb|> (basis (K 258) 4)|\empim
 \bmpi\verb|(<<AlLp[4 <<AlLp[5 6]>>]>> <<AlLp[2 <<AlLp[3 4]>>][2 <<AlLp[3 4]>>]>>)|\empi
 \bmpi\verb|> (length *)|\empim
 \bmpi\verb|2|\empix
 \bmpi\verb|> (basis (K 17) 4)|\empim
 \bmpi\verb|Error: attempt to call `:LOCALLY-EFFECTIVE' which is an undefined function.|\empi
 \emp

The basis in dimension 4 of the chain complex \#\boxtt{K258} is computed, it is
a list of length 2 (`\boxtt{*}' = the last returned object). The elements of
the basis themselves are ``algebraic loops'' (\boxtt{AlLp}), elements of some
appropriate cobar constructions.

On the contrary you see it is not possible to obtain the basis in dimension 4
of the chain complex \#\boxtt{K17} = \(C_\ast X\); the necessary functional
object is in fact the keyword \boxtt{:locally-effective} which generates an
error.

A homology group of \(X\) can be computed:

 \bmp
 \bmpi\verb|> (homology X 5)|\empim
 \bmpi\verb|Homology in dimension 5 :|\empi
 \bmpi\verb|Component Z/4Z|\empi
 \bmpi\verb|Component Z/2Z|\empi
 \bmpi\verb|Component Z|\empi
 \bmpi\verb|---done---|\empi
 \emp

\noindent which means \(H_5 X = \bZ_4 \oplus \bZ_2 \oplus \bZ\). It is in fact
the homology of \#\boxtt{K258}:

 \bmp
 \bmpi\verb|> (homology (K 258) 5)|\empim
 \bmpi\verb|Homology in dimension 5 :|\empi
 \bmpi\verb|Component Z/4Z|\empi
 \bmpi\verb|Component Z/2Z|\empi
 \bmpi\verb|Component Z|\empi
 \bmpi\verb|---done---|\empi
 \emp

The strong chain equivalence \#\boxtt{K268} contains three chain complexes
\emph{and} two reductions, therefore four chain complex morphisms and two
homotopy operators. In particular there is in the right reduction \(\rho_r:
K_{256} \rrdc K_{258}\) a \underline{r}ight \(g: K_{258} \rightarrow K_{256}\)
reachable by means of a \boxtt{\underline{r}g} function in the program; in the
same way the left reduction \(\rho_\ell: K_{17} \lrdc K_{256}\) contains a
\underline{l}eft \(f: K_{256} \rightarrow K_{17}\) reachable thanks to a
\boxtt{\underline{l}f} function. The Kenzo program can use these maps for
arbitrary generators or combinations. For example the next Lisp statements play
to verify the composite of the left \(f\) and the right \(g\) is compatible
with the differentials.

We assign to the symbol \boxtt{gen} the first generator of \#\boxtt{K258} in
dimension 4, we apply the right \(g\) (\boxtt{rg}) to this generator, then the
left \(f\) (\boxtt{lf}), finally the differential of \#\boxtt{K17}:

 \bmp
 \bmpi\verb|> (setf gen (first (basis (K 258) 4)))|\empim
 \bmpi\verb|<<AlLp[4 <<AlLp[5 6]>>]>>|\empix
 \bmpi\verb|> (rg SCE 4 gen)|\empim
 \bmpi\verb|----------------------------------------------------------------------{CMBN 4}|\empi
 \bmpi\verb| <-1 * <BcnB <TnPr <<AlLp[4 <<Loop[2-1 4][4-3 4]>>]>> <TnPr ... ...|\empi
 \bmpi\verb| <1 * <BcnB <TnPr ... ...|\empi
 \bmpi\verb|------------------------------------------------------------------------------|\empix
 \bmpi\verb|> (lf SCE *)|\empim
 \bmpi\verb|----------------------------------------------------------------------{CMBN 4}|\empi
 \bmpi\verb| <-2 * <<Loop[1-0 <<Loop[4]>>][3-2 <<Loop[4]>>]>>>|\empi
 \bmpi\verb| <2 * <<Loop[2-0 ...|\empi
 \bmpi\verb| [... Lines deleted...]|\empi
 \bmpi\verb|------------------------------------------------------------------------------|\empix
 \bmpi\verb|> (? (K 17) *)|\empim
 \bmpi\verb|----------------------------------------------------------------------{CMBN 3}|\empi
 \bmpi\verb|<-2 * <<Loop[<<Loop[3 4][5]>>]>>>|\empi
 \bmpi\verb|------------------------------------------------------------------------------|\empi
 \emp

The result is an actual ``loop of loops'', \((-2) \times\) some simplex in
\(X\). Large parts of the intermediate results are not showed. A result between
two dash lines `\boxtt{---}' labeled for example \boxtt{\{CMBN 3\}} is a
combination of degree 3 of integer coefficients and generators, one term per
line.

The other path consists in applying to the same generator firstly the
differential of \#\boxtt{K258} and then the same maps:

 \bmp
 \bmpi\verb|> (? (k 258) 4 gen)|\empim
 \bmpi\verb|----------------------------------------------------------------------{CMBN 3}|\empi
 \bmpi\verb|<2 * <<AlLp[3 <<AlLp[4 5]>>]>>>|\empi
 \bmpi\verb|------------------------------------------------------------------------------|\empix
 \bmpi\verb|> (rg sce *)|\empim
 \bmpi\verb|----------------------------------------------------------------------{CMBN 3}|\empi
 \bmpi\verb|<2 * <BcnB <TnPr <<AlLp[3 <<Loop[3 4][5]>>]>> <TnPr <<Loop>> <<Loop>>>>>>|\empi
 \bmpi\verb|<-2 * <BcnD <<AlLp[3 <BcnB <TnPr <<AlLp[4 5]>> <TnPr 0 <<Loop>>>>>]>>>>|\empi
 \bmpi\verb|<2 * <BcnD <<AlLp[3 <BcnD <<AlLp[4 5]>>>]>>>>|\empi
 \bmpi\verb|------------------------------------------------------------------------------|\empix
 \bmpi\verb|> (lf sce *)|\empim
 \bmpi\verb|----------------------------------------------------------------------{CMBN 3}|\empi
 \bmpi\verb|<-2 * <<Loop[<<Loop[3 4][5]>>]>>>|\empi
 \bmpi\verb|------------------------------------------------------------------------------|\empi
 \emp

The results are the same.

[Section to be continued, see Section~\ref{jmahz}]

\section{The fundamental theorem of Effective Homology.}

It is well known (?) the classical spectral sequences (Serre, Eilenberg-Moore,
Adams, \ldots) \emph{are not} algorithms. See for
example~\cite[Section~1.1]{MCCL}, in particular the comments following the
unique theorem of the quoted section: most often, the available input for a
spectral sequence \emph{does not} determine the higher differentials. Something
more is necessary for this essential problem, and it happens the category
\(\SSEH\) is from this point of view a perfect solution; moreover it is a
simple solution, once the possibilities of functional programming are
understood.

\begin{mth}\label{qxtmx} ---
Let \[\chi: (X_i)_{1 \leq i \leq n} \mapsto Y\] be a ``reasonable''
construction of the Algebraic Topology world producing \(Y\) from
the~\(X_i\)'s. Then an algorithm \(\chi_{EH}\) can be written down which is a
version with effective homology of the construction \(\chi\):
\[\chi_{EH}: ((X_i)_{EH})_{1 \leq i \leq n} \mapsto Y_{EH}\]
\end{mth}

Most often the \(X_i\)'s and \(Y\) are topological spaces. A construction is
``reasonable'' if it leads to some classical spectral sequence giving to the
topologists the illusion that if the homology (for example) of the \(X_i\)'s is
known, then the homology of \(Y\) can be ``deduced''.

A typical and important situation of this sort is the case where \(X\) is a
simply connected space and \(\chi = \Omega\) is the loop space functor: \(Y =
\Omega X\). The Eilenberg-Moore spectral sequence gives interesting relations
between \(H_\ast X\) and \(H_\ast \Omega X\), but this spectral sequence
\emph{is not} an algorithm computing \(H_\ast \Omega X\) from \(H_\ast X\), for
a simple reason: it is possible \(H_\ast X = H_\ast X'\) and \(H_\ast \Omega X
\neq H_\ast \Omega X'\). More precisely, the cobar construction~\cite{ADHL}
gives the homology of the \emph{first} loop space when some coproduct is
available around \(H_\ast X\), but the cobar construction \emph{does not} give
a coproduct around \(H_\ast \Omega X\), so that the process cannot be iterated;
this is \emph{Adams' problem}: how to iterate the cobar construction? More than
twenty years after Adams, Baues succeeded in a beautiful work~\cite{BAUS1} in
iterating \emph{one time} the cobar construction, giving the homology of the
second loop space \(\Omega^2 X\) in reasonable situations, but Baues' method
cannot be extended for the homology of \(\Omega^3 X\) either.

The category \(\SSEH\) gives at once a complete and simple solution for Adams'
problem; it is a consequence of the following particular case of
Meta-Theorem~\ref{qxtmx}.

\begin{thr}\label{jovfr} ---
An algorithm \(\Omega_{EH}\) can be written down:
\[
\Omega_{EH} : X_{EH} \mapsto (\Omega X)_{EH}
\]
producing a version \emph{with effective homology} of the loop space \(\Omega
X\) when a version \emph{with effective homology} of the initial simply
connected space \(X\) is given.
\end{thr}

The algorithm \(\Omega_{EH}\) not only \emph{can be} written down, but \emph{it
is} written down; the algorithm \(\Omega_{EH}\) is certainly the most important
component of the Kenzo program~\cite{DRSS}, a program which is by itself the
most detailed proof which can be required for Theorem~\ref{jovfr}, of course
not very convenient for an ordinary reader\footnote{Several articles containing
such a proof written in common mathematical language have been proposed to
various mathematical journals, but they were always rejected by the editorial
boards, see in particular~\cite{SRGR6}. It is likely that the totally new
nature of the result, which can be stated and proved only in a computational
framework, does not fit the standard style expected by the referees. A matter
of evolution; yet the scientific journals should precisely be mainly interested
by the papers reflecting new unavoidable scientific trends, papers which of
course are a little more difficult to be appropriately
refereed.\\See~\cite{SRGR4} for a survey which gives the plan and the main
ideas of the proof of Theorem~\ref{jovfr}.}.

The data type of the output \((\Omega X)_{EH}\) is exactly the same as the data
type of the input \(X_{EH}\), so that the algorithm \(\Omega_{EH}\) can be
trivially iterated.

\begin{thr}\textbf{\emph{(Solution of Adams' problem\footnote{This 
theorem \emph{is} also a solution for Carlsson's and Milgram's problem~\cite[p.
545, Section~6]{CRML}, a problem the authors cannot properly state, again
because of the lack of a computational framework.})}}
--- An algorithm \textbf{\emph{ICB}} (\underline{i}terated
\underline{c}o\underline{b}ar) can be written down:
\[
\textbf{\emph{ICB}}: (X_{EH}, n) \mapsto (\Omega^n X)_{EH}
\]
which produces a version \emph{with effective homology} of the \(n\)-th loop
space \(\Omega^n X\) when a version \emph{with effective homology} of the
initial space \(X\), assumed to be \(n\)-connected, is given.
\end{thr}

When \(X_{EH} = (X, C_\ast X, EC_\ast^X, \varepsilon^X)\) is given, the
algorithm \textbf{ICB} produces a 4-tuple \((\Omega^n X)_{EH} = (\Omega^n X,
C_\ast \Omega^n X, EC_\ast^{\Omega^n X}, \varepsilon^{\Omega^n X})\), where the
``\(n\)-th cobar'' of \(EC_\ast^X\) is the third component \(EC_\ast^{\Omega^n
X}\). This \(n\)-th cobar cannot be constructed from \(EC_\ast^X\) only; the
first cobar needs the coproduct of \(C_\ast X\) and the \(n\)-th cobar needs
much more supplementary informations hidden in \(X\) and \(\varepsilon^X\);
these objects \(X\) and \(\varepsilon^X\) are locally effective and model
mathematical objects which are infinite; yet \(X\) and \(\varepsilon^X\) are
\emph{finite} machine objects (pleonasm), namely finite bit strings actually
created, processed and used by the Kenzo program; this process works thanks to
\emph{functional programming}.

The further components \(\Omega^n X\) and \(\varepsilon^{\Omega^n X}\) in the
result would allow to undertake other calculations starting from \(\Omega^n
X\).

\section{A small machine demonstration [sequel].}\label{jmahz}

Let us consider again the space \(X\) of Section~\ref{czlgy}. The Kenzo program
had constructed a version \emph{with effective homology} of this space,
allowing in particular to compute its homology groups. Much more important,
because of Theorem~\ref{jovfr}, the machine object \(\Omega_{EH}\) of the same
program can be applied to produce a version \emph{with effective homology} of
the loop space \(\Omega X\):

 \bmp
 \bmpi\verb|> (setf OX (loop-space X))|\empim
 \bmpi\verb|[K273 Simplicial-Group]|\empi
 \emp

Kenzo\footnote{The Kenzo function \boxtt{loop-space} follows the modern rules
of \emph{Object Oriented Programming}~(OOP): if the argument is a simplicial
set, then the Kan model of the loop space is constructed, and if furthermore
the argument contains the effective homology of the initial simplicial set,
then the \boxtt{loop-space} function constructs also the effective homology of
the loop space.} returns a new locally effective simplicial \emph{group}, the
Kan model of~\(\Omega X\) and its effective homology:

 \bmp
 \bmpi\verb|> (setf SCE2 (efhm OX))|\empim
 \bmpi\verb|[K405 Equivalence K273 <= K395 => K391]|\empi
 \emp

\noindent with which exactly the same experiences which were tried with the
effective homology of \(X\) in Section~\ref{czlgy} could be repeated. In
particular the right chain complex \#\boxtt{K391} is effective and allows a
user to compute a homology group:

 \bmp
 \bmpi\verb|> (homology OX 5)|\empim
 \bmpi\verb|Component Z16|\empi
 \bmpi\verb|Component Z8|\empi
 \bmpi\verb|...|\empi
 \bmpi\verb|...|\empi
 \bmpi\verb|Component Z2|\empi
 \bmpi\verb|Component Z2|\empi
 \emp

\noindent where we have deleted 21 lines, for the result is in fact:
\[
H_5(\Omega(\Omega(\Omega (P^\infty(\bR)/P^3(\bR)) \cup_4 D^4) \cup_2 D^3)) =
\bZ_{16} \oplus \bZ_8 \oplus \bZ_2^{23}.
\]

This is nothing but the corresponding homology group of \#\boxtt{K391}.

To our knowledge, the Kenzo program is the only object, human or not, currently
able to reach this result. See also~\cite{SCHN,SMTH1} for two other interesting
theoretical solutions which unfortunately have not yet led to concrete machine
programs.

Adams' and Carlsson-Milgram's problems are solved.

\section{The category \(SSEH\) and the Postnikov category.}\label{deoad}

A rough ``definition'' of the Postnikov category was given in
Section~\ref{mssfy}, but we must be now more precise to obtain a correct
relation between the category \(\SSEH\) and the Postnikov category.

\begin{dfn}
\emph{--- An \emph{Abelian group of finite type} \(\pi\) is a direct sum \(\pi
= \bZ/d_1 \bZ \oplus \cdots \oplus \bZ/d_n\bZ\) where every \(d_i\) is a
non-negative integer and \(d_{i-1}\) divides \(d_i\). We denote by \(\Pi\) the
\emph{set} of \emph{these} groups}.
\end{dfn}

The set \(\Pi\) is designed for having exactly \emph{one} group isomorphic to
every Abelian group of finite type. For example the group \(H_5 X\) in
Section~\ref{czlgy} is \emph{isomorphic} to the element of \(\Pi\) defined by
the integer triple \((2,4,0)\), that is, the group \(\bZ_2 \oplus \bZ_4 \oplus
\bZ\), but there are 128 different isomorphisms.

\begin{dfn}
\emph{--- A \emph{Postnikov tower} is a pair of sequences \(((\pi_n)_{n \geq
2}, ( k_n)_{n \geq 3})\) where \(\pi_n \in \Pi\) and \(k_n \in H^{n+1}(X_{n-1},
\pi_n)\), if \(X_n\) is the \(n\)-th stage of the Postnikov tower constructed
by the standard process.}
\end{dfn}

Because \(\pi_n\) \emph{is} some \emph{precise} group, the standard
Eilenberg-MacLane process gives a \emph{precise} \(K(\pi_n, n)\), a Kan
simplicial set, producing in turn a \emph{precise} \(n\)-th Postnikov stage
\(X_n\) and with the next \(\pi_{n+1}\) a \emph{precise} cohomology group
\(H^{n+2}(X_n, \pi_{n+1})\) \emph{where} the \(k_{n+1}\) \emph{must} be
``chosen''. A Postnikov tower so produces in a deterministic way a realization
\(X\). A morphism \(f: (\pi_n, k_n) \rightarrow (\pi'_n, k'_n)\) between two
Postnikov towers is a collection \((f_n: \pi_n \rightarrow \pi'_n)\) of group
homomorphisms compatible with the \(k_n\)'s and \(k'_n\)'s, that is, satisfying
Baues' relation, and we have so defined the Postnikov category \(\cP\). The
isomorphism problem consists in deciding whether two Postnikov towers \((\pi_n,
k_n)\) and \((\pi'_n, k'_n)\) produce realizations with the same homotopy type,
that is, because of the context, that are isomorphic. Of course the condition
\(\pi_n = \pi'_n\) is required for every~\(n\), but simple examples show that
the condition \(k_n = k'_n\) on the contrary is not necessarily required. This
is the reason why the \(k_n\)'s \emph{are not} invariants of the homotopy type.

The computable category \(\SSEH\) makes the realization process computable.

\begin{thr}\label{lkdyu}
--- An algorithm \textbf{\emph{PR}} (Postnikov realization) can be written down:
\[
\textbf{\emph{PR}}: \cP \rightarrow \SSEH
\]
implementing the realization process.
\end{thr}

In fact the situation is significantly more complex. Before being\ldots
\emph{true}, the statement of this theorem must \emph{make sense}, so that a
machine implementation of the category \(\cP\) must be available, at least from
a theoretical point of view. This is obtained thanks to the category \(\SSEH\)
itself: a component \(k_n\) must be a machine object, so that the \emph{data
type} \(H^{n+1}(X_{n-1}, \pi_n)\) where \(k_n\) is to be picked up must be
\emph{previously} defined, which is possible only if a calculation of this
cohomology group can be undertaken. And again it is the category \(\SSEH\)
which gives this possibility. It is an amusing situation where a category, the
category \(\SSEH\), is simultanously used to give sense to the \emph{statement}
of a theorem, and \emph{synchronously} finally to prove it.

Combining Theorem~\ref{lkdyu} with the appropriate particular cases of
Meta-Theorem~\ref{qxtmx}, we see the problem implicitly stated in the framed
title inscription is now solved. In particular, the Kenzo program allows a
Postnikov user to undertake many computations of this sort.

It is not possible, with the currently available tools, to make the categories
\(\cP\) and \(\SSEH\) \emph{effectively} equivalent, of course up to the
homotopy relation.

\begin{thr}\label{rohdy}
--- An algorithm \textbf{\emph{SP}} can be written down:
\[
 \textbf{\emph{SP}}:\ \ \SSEH \rightarrow \cP \times \cI\ \ :\ \
 X \mapsto (\pi_n, k_n, I_n)_{n \geq 2}
\]
where the component \(I_n\) is \emph{some} isomorphism \(I_n: \pi_n(X) \cong
\pi_n \in \Pi =\) the set of ``canonical'' models of Abelian groups of finite
type; the \(k\)-invariants \(k_n\) are unambiguously defined \emph{only when}
the \(I_n\)'s are chosen.
\end{thr}

The algorithm \textbf{SP} is essentially non-unique, for the \emph{choice} of
the component~\(I_n\) is arbitrary, and the various choices of these
isomorphisms will produce all the possible collections of \(k\)-``invariants''
(\(!\)). Unfortunately the group \(GL(p,\bZ)\) for example is infinite for \(p
> 1\), so that it is a non-trivial arithmetical problem to determine whether
two collections of \(k\)-invariants correspond to the same homotopy type or
not. For example the problem has an obvious solution up to arbitrary dimensions
if every~\(\pi_n\) is finite, but as soon as a \(\pi_n\) is not finite, we are
in front of interesting but difficult problems of arithmetic, to our knowledge
not yet solved in general\footnote{Compare with~\cite[pp. 54-59]{SCHN}; the
possible equivalence of \(k_n\) and \(k'_n\) with respect to some automorphism
of the \emph{last} \(\pi_n\) is there proved decidable, which is relatively
easy. But this does not seem to be sufficient, because the possible
automorphisms of \emph{all} the previous \(\pi_m\), \(m \leq n\), must be
considered. The example of a Postnikov tower where only \(\pi_2 = \bZ^p\) and
\(\pi_5 = \bZ\) given in Section~\ref{mssfy} shows the main problem for the
equivalence of \(k\)-invariants is in the automorphisms of \(\pi_2\), because
the group of automorphisms is \(GL(p,\bZ)\), leading to essential hard
arithmetical problems. In a later preprint, not published, Rolf Sch\"on
considers again the problem, solves it in the case where all the \(\pi_n\)'s
are finite, and \emph{announces} a general solution which ``takes considerable
work''. The authors have not succeeded in getting in touch with Rolf Sch\"on
for several years, and any indication about his current location would be
welcome. Note that these comments in particular cancel the assertion
in~\cite[Section 5.4]{SRGR5} about the classification problem, which assumed
the correctness of Sch\"on's paper: the solutions called JS, SRH and SRG
in~\cite{SRGR5} solve only the \emph{soft problem}; an essential gap remains
present for the \emph{hard problem}.}.

In conclusion, thanks to the computable category \(\SSEH\), the category
\(\cP\) \emph{becomes} also a computable category. There are ``good'' but
non-canonical correspondances between these categories. Both categories solve
the \emph{soft problem} and from this point of view give equivalent results.
Both categories would solve the \emph{hard problem} if the equivalence problem
between systems of \(k\)-``invariants'' was \emph{effectively} solved.

It is easier now to understand the common confusion about the nature of the
\(k\)-invariants. We follow exactly here \cite[\S25]{MAY}, up to obvious slight
differences of notations. If \(X\) is a topological space, we can start with a
minimal Kan model of \(X\), ``unique'' up to numerous \emph{different}
isomorphisms in general; the Postnikov stages \(X_{n-1}\) and \(X_n\) are then
\emph{canonical} quotients of~\(X\). There is also a \emph{canonical} fibration
between \(X_n\) and \(X_{n-1}\), the fiber space of which being the
\emph{canonical} space \(K(\pi_n\underline{(X)}, n)\), defining
\emph{unambiguously} a \[k_n \in H^{n+1}(X_{n-1}, \pi_n\underline{(X)}) =
H^{n+1}(\underline{X}/\sim_{n-1}, \pi_n\underline{(X)}).\] It is then clear
that a claimed invariant living in \(H^{n+1}(X_{n-1}, \pi_n)\) with \(\pi_n \in
\Pi\) depends on an isomorphism \(\pi_n \cong \pi_n(X)\), which is essentially
the \(g\) correctly considered at~\cite[Theorem 25.7]{MAY}. But if you think
that \(X_{n-1}\) comes \emph{only} from the previous data \(\pi_2, \pi_3, k_3,
\ldots, \pi_{n-1}, k_{n-1}\) and not from~\(X\) itself, you can \emph{freely}
apply a self-equivalence of \(X_{n-1}\) to change~(\(!\)) the invariant, the
\emph{fibration} \(X_n \rightarrow X_{n-1}\) is changed, but on the contrary
the homotopy type of \(X_n\) remains unchanged: \emph{different} invariants
correspond to \emph{equal homotopy types}. The group of all the
self-equivalences of \(X_{n-1}\) must be considered, of course in general a
serious question. The only way to cancel this ambiguity consists in choosing a
well defined partial equivalence between \(X\) and \(X_{n-1}\), which amounts
to choosing some isomorphisms \(\pi_i \cong \pi_i(X)\) for \(2 \leq i < n\).

Maybe it is useful to recall an \emph{invariant} with respect to any notion
must be chosen in a ``fixed'' world independent of the object the invariant of
which is being defined. Otherwise the definitively simplest complete invariant
for the homotopy type of \(X\) is \(X\) itself. Not very interesting. The
\emph{non-ambiguous} definition of \(k_n\) above lives in a set the definition
of which contains two occurences of \(X\) and this is forbidden when an
\emph{invariant} of \(X\) is defined. And which must be called an error
in~\cite[Theorem~25.7]{MAY} comes from the definition \(k_n \in
H^{n+1}(X_{n-1}, \pi_n)\) (p.~113, line~19), and the notation \(\pi_n =
\pi_n(X, \emptyset)\) (line~14), again two illegal occurences of \(X\) in this
situation. It is the reason why, in the statement of Theorem~\ref{rohdy},
\(\pi_n\) is not \emph{equal} to \(\pi_n(X)\), they are only isomorphic through
some isomorphism which plays an essential role.

The classical example of the minimal polynomial of a matrix is helpful; if the
ground field \(K\) is given, then the minimal polynomial can be chosen once and
for all in~\(K[\lambda]\), a set of polynomials independent of the particular
considered matrix. So that if two matrices are conjugate, more generally if two
endomorphisms of two finite-dimensional \(K\)-vector spaces are
\emph{conjugate}, their minimal polynomials are \emph{equal}, not mysteriously
``\emph{isomorphic}''; this is the reason why the miminal polynomial is a
correct conjugation \emph{invariant}.

\subsection{Localization.}\label{nnmrf}

Another natural idea must also be considered. It is usual to split a
topological problem~\(P\) into a rational problem \(P_0\) and \(p\)-problems
\((P_p)_{p \in \mathfrak{P}}\) for \(p\) running the prime numbers
\(\mathfrak{P}\). Then a solution for every problem \(P_p\)
\emph{theoretically} produces a solution for the initial problem \(P\). Let us
take again our fetish example of the Postnikov towers with only \(\pi_2 =
\bZ^k\) and \(\pi_{2n-1} = \bZ\) where the \(k_{2n-1}\)-invariant is \(\chi \in
S_k^n\), the \(\bZ\)-module of the homogeneous \(\bZ\)-polynomials of
degree~\(n\) with respect to \(k\) variables. The localization method produces
a localized Postnikov tower \(T_p\) for every element \(p \in \{0\} \cup
\mathfrak{P}\). We are then in front of a list of problems.

\begin{prb}\label{pczfl} ---  Let \(\chi\) and \(\chi'\) be  two
polynomials defining Postnikov towers \(T\) and~\(T'\), producing in turn
families of Postnikov towers \((T_p)\) and \((T'_p)\).
\begin{enumerate}
\item
Let \(p\) be an element of \(\{0\} \cup \mathfrak{P}\); can the mapping \(T
\mapsto T_p\) be made \emph{effective}?
\item
Can the isomorphism problem between \(T_p\) and \(T'_p\) be \emph{effectively}
solved?
\item
Does there exist an integer \(p_0\) and an argument allowing to exempt us from
this study for \(p > p_0\)?
\item
If the isomorphism problem between \(T_p\) and \(T'_p\) has a positive
effective solution for every \(p\), does there exist an \emph{effective}
process allowing to construct an isomorphism between \(T\) and \(T'\)?
\end{enumerate}
\end{prb}

Subproblems 1 and 2 of Problem~\ref{pczfl} probably are ``exercises'', but the
subproblems~3 and~4 seem serious. Notice again a general solution for
Problem~\ref{pczfl} would solve at once the problem of the \(\bZ\)-linear
equivalence between elements of \(S_k^n\) for every \(n\), while the
arithmeticians currently know the solution only for \(n = 2\), a non-trivial
problem~\cite{WTSN}. From this point of view, it would be interesting to
\emph{translate} the known solution which is available for \(n = 2\) into a
solution of Problem~\ref{pczfl}.

\section{The category \(\SSEH\) and the \(E_\infty\)-algebras.}\label{bxwta}

We had briefly mentioned in Section~\ref{bdpci} the possibility of other
solutions for the \emph{hard problem} based on \(E_\infty\)-operads.

A particularly interesting \(E_\infty\)-operad is the \emph{surjection operad}
\(\cS\) defined and studied in~\cite{BRFR}, a work undertaken to make
completely explicit\footnote{The considerations of Section~\ref{nnmrf} can
again be applied when comparing the Steenrod operad~\(\cS\), a \(\bZ\)-operad,
and the ``abstract'' \(p\)-localized operads of Mandell.} some results of
Mandell's paper~\cite{MNDL} already quoted in Section~\ref{bdpci}. The
so-called surjection operad and its action on a simplicial set can be
understood as a ``complete'' generalization of the Steenrod operations, and we
therefore propose to call it the \emph{Steenrod operad}, which furthermore
allows to naturally keep the same notation \(\cS\).

\begin{thr}\label{lvuop}
--- A functorial algorithm \textbf{\emph{SSC}} (simplicial sets to Steenrod chain complexes)
can be written down:
\[
\textbf{\emph{SSC}}: \SSEH \rightarrow \CCS
\]
where \(\CCS\) is the category of the free \(\bZ\)-chain complexes of finite
type provided with a \(CB\cS\)-operadic structure.
\end{thr}

An appropriate bar construction can be applied to the operad \(\cS\) to produce
a cooperad \(B\cS\); then an analogous cobar construction can in turn be
applied to this cooperad to produce a new operad denoted by \(CB\cS\), another
model for an \(E_\infty\)-operad which has the following advantage\footnote{We
would like to thank Tornike Kadeishvili for his clear and useful explanations
about this process.}: let \(f: C_\ast \rightarrow D_\ast\) be a chain
equivalence between two free \(\bZ\)-chain complexes; then every
\(CB\cS\)-structure on \(C_\ast\) induces such a structure on \(D_\ast\).

\begin{dfn}
\emph{--- A \emph{Steenrod chain complex} is a free \(\bZ\)-chain complex
provided with a \(\cS\)-structure or with a \(CB\cS\)-structure}.
\end{dfn}

Let \(X\) be an object of \(\SSEH\), that is, a simplicial set with effective
homology. The article~\cite{BRFR} explains how the initial definition by
Steenrod of his famous cohomological operations can be naturally used to
install a canonical \(\cS\)-structure on the chain complex \(C_\ast X\); the
strong chain equivalence \(\varepsilon^X: C_\ast X \eqvl EC_\ast^X\) then
allows to install a \(CB\cS\)-structure on \(EC_\ast^X\), and this is enough to
define the functor \textbf{SSC}.

Taking account of Mandell's article~\cite{MNDL}, the following problems are
natural.

\begin{prb}\label{cnzfs} ---
Does there exist an \emph{algorithm} \textbf{\emph{R}} (realizability):
\[
\textbf{\emph{R}}: \CCS \times \bN \rightarrow \textbf{\emph{Bool}} =
\{\top,\bot\}
\]
allowing to decide whether some object \(C_\ast \in \CCS\) corresponds or not
to some topological object up to some given dimension?
\end{prb}

Because of the Characterization Theorem~\cite[p. 2]{MNDL}, a solution for this
problem is probably a ``simple'' exercise, simple in theory but the operad
\(CB\cS\) is rather sophisticated, so that a \emph{concrete} solution seems a
nice challenge. Furthermore the Characterization Theorem is stated and proved
in characteristic~\(p\) and obtaining the analogous result with respect to the
ground ring~\(\bZ\) could be a little difficult.

\begin{prb}\label{qxzqq} ---
Does there exist an \emph{algorithm} \textbf{\emph{SHT}} (same homotopy type):
\[
\textbf{\emph{SHT}}: \CCS' \times \CCS' \times \bN \rightarrow
\textbf{\emph{Bool}}
\]
allowing to decide whether two \(\CCS\)-objects obtained through the
\textbf{\emph{SSC}}-algorithms, therefore certainly corresponding to actual
recursive simplicial sets, have the same homotopy type or not, of course up to
some given dimension?
\end{prb}

The authors are not sufficiently experienced in operadic techniques to estimate
the difficulty of this question. The Main Theorem of~\cite[p. 1]{MNDL} seems to
imply that the same considerations as for Problem~\ref{cnzfs} could be applied;
but as already observed, an effective solution of Problem~\ref{qxzqq} would
indirectly solve crucial computability problems in arithmetic, problems which
seem to raise essential obstacles in front of the professionals. It is
difficult to think the \(E_\infty\)-operad could be a mandatory tool to solve
these arithmetical problems, so that for a concrete solution it is more
tempting to solve directly the arithmetical problems and to use \emph{only} the
category \(\SSEH\) and its satellite category \(\cP\), thanks to
Theorems~\ref{lkdyu} and~\ref{rohdy}, to obtain a solution of the hard problem.

If the operadic methods become unavoidable, it seems terribly difficult to
design \emph{directly} the category \(\CCS\) as a \emph{computable} category.
We think it would be more sensible to work \emph{simultaneously} with the
categories \(\SSEH\) and \(\CCS\): it is frequent in mathematics in general,
and in computer science in particular, that it is not a good idea to give up
too early informations which \emph{look} redundant. This is well known for
example by the theoreticians in homotopy theory: it is much better to work with
an \emph{explicit} homotopy equivalence than only with the \emph{existence} of
such an object, and it is still better to keep also the various maps which
describe how this homotopy equivalence actually is one, and so on. This is
nothing but the philosophy always underlying when we work with
\(E_\infty\)-operads.

In our situation, Theorem~\ref{lvuop} implies a simplicial set with effective
homology \(X_{EH}\) \emph{contains} in an effective way a Steenrod chain
complex; and we do not need any realizability criterion, an object of \(\SSEH\)
\emph{certainly} corrresponds to a \emph{genuine} topological space. Therefore
the right objects to work with in Algebraic Topology could be the pairs
\((X_{EH}, \Sigma_\cS^X)\) where the second component \(\Sigma_\cS^X\) is the
\(CB\cS\)-structure induced on \(EC_\ast^X\) by the canonical Steenrod
structure on \(C_\ast X\). Then, when a new object is constructed from such
objects, the ingredients present in the second components could facilitate the
computation of some parts of the constructed object, but others would certainly
be obtained much more easily thanks to the first components.

In a sense the success of the category~\(\SSEH\) is already of this sort:
instead of working only with a chain complex \(EC_\ast^X\) describing the
homology of \(X\), certainly in general non-sufficient for the planned
computations, it is much better to work with~\(X\) itself under its locally
effective form, the only form which can be processed on a machine when \(X\) is
not of finite type. The amazing fact is that this is sufficient to solve many
computability problems, though this version of \(X\) \emph{does not
effectively} defines the mathematical object \(X\), because of G\"odel and his
friends, see~\cite[Section~5.3]{SRGR5}. The same people, helped by
Matiyasevich~\cite{MTSV}, have also made impossible a universal solver of
systems of polynomial \(\bZ\)-equations, and after all, the hard problem is
equivalent to a problem about such equations, so that we cannot even be sure,
up to further information, a solution of the hard problem~\emph{exists}.

\vspace{20pt}

 \noindent{\footnotesize{\begin{tabular}{rl}
 JR: & \texttt{Julio.Rubio@dmc.unirioja.es} \\[5pt]
 FS: & \texttt{Francis.Sergeraert@ujf-grenoble.fr}
\end{tabular}}}


\begin{thebibliography}{99}
\bibitem{ADHL}
   J. F. Adams, Peter J. Hilton.
   \emph{On the chain algebra of a loop space}.
   Commentarii Mathematici Helvetici, 1956, vol. 30, pp 305-330.
\bibitem{BAUS1}
   Hans J. Baues.
   \emph{Geometry of loop spaces and the cobar construction}.
   Memoirs of the American Mathematical Society,
   1980, vol. 230.
\bibitem{BAUS4}
   Hans J. Baues.
   \emph{Homotopy types}.
   in \cite{JAMS}, pp 1-72.
\bibitem{BRFR}
   Clemens Berger and Benoit Fresse.
   \emph{Combinatorial operad actions on cochains}.
   Preprint.
\bibitem{CRML}
   Gunnar Carlsson and R. James Milgram.
   \emph{Stable homotopy and iterated loop spaces}.
   in \cite{JAMS}, pp 505-583.
\bibitem{DRSS}
   Xavier Dousson, Julio Rubio, Francis Sergeraert and Yvon Siret.
   \emph{The Kenzo program}.
   \texttt{http://www-fourier.ujf-grenoble.fr/\~{ }sergerar/Kenzo/}
\bibitem{ELMC1}
   Samuel Eilenberg, Saunders MacLane.
   \emph{On the groups $H(\pi,n)$, I}.
   Annals of Mathematics, 1953, vol. 58, pp 55-106.
\bibitem{ELMC2}
   Samuel Eilenberg, Saunders MacLane.
   \emph{On the groups $H(\pi,n)$, II}.
   Annals of Mathematics, 1954, vol. 60, pp 49-139.
\bibitem{EDM}
   \emph{Encyclopedic Dictionary of Mathematics}, sub-article \emph{Postnikov
   complexes}, in different articles according to the edition (look for
   \emph{Postnikov complex} in the final index). Mathematical Society of Japan and
   American Mathematical Society.
\bibitem{JAMS}
   \emph{Handbook of Algebraic Topology} (Edited by I.M. James).
   North-Holland (1995).
\bibitem{MNDL}
   M. Mandell.
   \emph{\(E_\infty\)-algebras and \(p\)-adic homotopy theory}.
   Topology, 2001, vol. 40, pp 43-94.
\bibitem{MTSV}
   Yuri Matiyasevich.
   \emph{Hilbert's tenth problem}.
   MIT Press, 1993.
\bibitem{MAY}
   J. Peter May.
   \emph{Simplicial objects in algebraic topology}.
   Van Nostrand, 1967.
\bibitem{MCCL}
   John McCleary.
   \emph{User's guide to spectral sequences.}
   Publish or Perish, Wilmington DE, 1985.
\bibitem{RBSR6}
   Julio Rubio, Francis Sergeraert.
   \emph{Constructive Algebraic Topology}.
   Bulletin des Sciences Math\'ematiques, 2002, vol. 126, pp 389-412.
\bibitem{RDBR}
   David E. Rydeheard, Rod M. Burstall.
   \emph{Computational Category Theory}.
   Prentice Hall, 1988.
\bibitem{SCHN}
   Rolf Sch\"on.
   \emph{Effective algebraic topology}.
   Memoirs of the American Mathematical Society,
   1991, vol. 451.
\bibitem{SRGR3}
   Francis Sergeraert.
   \emph{The computability problem in algebraic topology}.
   Advances in Mathematics, 1994, vol. 104, pp 1-29.
\bibitem{SRGR4}
   Francis Sergeraert.
   \emph{Effective homology, a survey}.\\
   \texttt{www-fourier.ujf-grenoble.fr/\~{ }sergerar/Papers/survey.dvi
   \mbox{\small or} ps}.
\bibitem{SRGR5}
   Francis Sergeraert.
   \emph{\(\maltese_K\), objet du 3\raisebox{.8ex}{\hspace{.1ex}\footnotesize e} type}.
   Gazette des Math\'ematiciens, 2000, vol. 86, pp 29-45.
\bibitem{SRGR6}
   Francis Sergeraert.\\
   \texttt{www-fourier.ujf-grenoble.fr/\~{ }sergerar/Papers/}.
\bibitem{SMRN}
   Vladimir A. Smirnov.
   \emph{The homology of iterated loop spaces}.
   Forum Mathematicum, 2002, vol.~14, pp ~345-381.
\bibitem{SMTH1}
   Justin R. Smith.
   \emph{Iterating the cobar construction}.
   Memoirs of the American Mathematical Society,
   1994, vol. 524.
\bibitem{TRVN}
   A.S. Troelstra, D. van Dalen.
   \emph{Constructivism in mathematics, an introduction}.
   North-Holland, 1988.
\bibitem{WTSN}
   G.L. Watson.
   \emph{Integral quadratic forms}.
   Cambridge University Press, 1960.
\end{thebibliography}
\end{document}